\documentclass[12pt,twoside]{article}

\begin{document}

\date{}




\centerline {\Large{\bf Related Fixed Point Theorems For Two Pairs }}

\centerline{}

\centerline {\Large{\bf  Of Mappings In Fuzzy Metric Spaces}}
\centerline {\Large{\bf }}

\centerline{\bf { T. K. Samanta , Sumit Mohinta and Iqbal H. Jebril }}
\centerline{}
\centerline{Department of Mathematics, Uluberia College,}
\centerline{Uluberia, Howrah, West Bengal, India-711315.} %
\centerline{e-mail: mumpu$_{-}$tapas5@yahoo.co.in}
\centerline{e-mail: sumit.mohinta@yahoo.com}
\centerline{Department of Mathematics, Taibah University, Saudi Arabia}
\centerline{e-mail: iqbal501@hotmail.com}
\begin{abstract}
\textbf{\emph{In this paper we establish the existence of related fixed points theorems for two pairs of mappings with different contraction conditions in two fuzzy metric spaces.}}
\end{abstract}


\newtheorem{Theorem}{\quad Theorem}[section]

\newtheorem{Definition}[Theorem]{\quad Definition}

\newtheorem{Corollary}[Theorem]{\quad Corollary}

\newtheorem{Lemma}[Theorem]{\quad Lemma}

\newtheorem{Note}[Theorem]{\quad Note}

\newtheorem{Remark}[Theorem]{\quad Remark}

\newtheorem{Result}[Theorem]{\quad Result}

\newtheorem{Proposition}[Theorem]{\quad Proposition}

\newtheorem{Example}[Theorem]{\quad Example}


\section{Introduction}

The concepts of fuzzy sets was initially investigated by Zadeh \cite{zadeh} in 1965 as a new way to represent imprecise facts or uncertainties or vagueness in everyday life. Subsequently , it was developed extensively by many authors and used in population dynamics , chaos control , computer programming , medicine , etc.
\\In 1975 Kramosil and Michalek \cite{Kramosil} introduced the concept of fuzzy metric spaces ( briefly , FM-spaces ) , which opened a new avenue for further development of analysis in such spaces. Later on it is modified and a few concepts of mathematical analysis have been developed by George and Veeramani \cite{Veeramani}. \\
Fisher \cite{Brian}, Aliouche and Fisher \cite{Aliouche}, Telci \cite{Telci} proved some related fixed point theorems in compact metric spaces. Recently, Rao et.al \cite{Rao} and \cite{Rao1} proved some related fixed point theorems in sequentially compact fuzzy metric spaces. However , the study of related fixed points for two pairs of mappings is also interesting. In this paper we extend this concept to fuzzy metric space and establish the existence of related fixed points theorems for two pairs of mappings.
This research modifies and generalizes the results of Fisher \cite{Brian},  R.K.Namdeo, S. Jain and B.Fisher \cite{Namdeo} under a different contraction
condition in two fuzzy metric spaces .

\section{Preliminaries}

We quote some definitions and statements of a few theorems which will be
needed in the sequel.

\begin{Definition}\cite{Schweizer}
A binary operation \, $\ast \; : \; [\,0 \; , \; 1\,] \; \times \;
[\,0 \; , \; 1\,] \;\, \longrightarrow \;\, [\,0 \; , \; 1\,]$ \, is
continuous \, $t$ - norm if \,$\ast$\, satisfies the
following conditions \, $:$ \\
$(\,i\,)$ \hspace{0.5cm} $\ast$ \, is commutative and associative ,
\\ $(\,ii\,)$ \hspace{0.3cm} $\ast$ \, is continuous , \\
$(\,iii\,)$ \hspace{0.2cm} $a \;\ast\;1 \;\,=\;\, a \hspace{1.2cm}
\forall \;\; a \;\; \varepsilon \;\; [\,0 \;,\; 1\,]$ , \\
$(\,iv\,)$ \hspace{0.3cm} $a \;\ast\; b \;\, \leq \;\, c \;\ast\; d$
\, whenever \, $a \;\leq\; c$  ,  $b \;\leq\; d$  and  $a \, , \, b
\, , \, c \, , \, d \;\, \varepsilon \;\;[\,0 \;,\; 1\,]$.
\end{Definition}

\begin{Definition}\cite{Veeramani}
The $3$-tuple $(\,X \,,\, \mu \,,\, \ast\,)$ is called a fuzzy metric space if X is an arbitrary non-empty set, \,$\ast$\, is a continuous t-norm and \,$\mu$ is a fuzzy set in \,$X^{2}\,\times\,(0,\infty)$\, satisfying the following conditions \,:
\\$(\,i\,)$ \hspace{0.5cm}$\mu\,(\,x \;,\,y \;,\; t\,) \;\,>\;\, 0 \, ;$ \\
$(\,ii\,)$ \hspace{0.2cm} $\mu\,(\,x \;,\,y\;,\; t\,) \;\,=\;\, 1$ \, if
and only if \, $x \;=\;y \,$\,
\\$(\,iii\,)$ \hspace{0.1cm}$\mu\,(\,x \;,\,y\;,\; t\,) \;\,=\;\,\mu\,(\,y \;,\,x\;,\; t\,);$
\\ $(\,iv\,)$\hspace{0.1cm} $\mu\,(\,x \;,\,y \;,\; s\,) \;\ast\; \mu\,(\,y \;,\,z \;,\; t\,)
\;\,\leq\;\, \mu\,(\,x \;,\,z \;,\; s\;+\;t\,
\,) \, ;$
\\$(\,v\,)$ \hspace{0.2cm}$\mu\,(\,x \,,\,y\,,\,\cdot\,) :(0 \,,\;\infty\,)\,\rightarrow
\,(0 \,,\;1]$ \, is continuous; \\ for all x\,,\,y\,,\,z  \,$\,\in\,X$\, and \,\,$t,\,s\,>\,0.$
\end{Definition}

\begin{Definition}\cite{Samanta}
Let \,$(\,X \,,\, \mu \,,\, \ast\,)$\, be a fuzzy metric space. A sequence \,$\{\,x_{\,n}\,\}_{\,n}$\, in \,$X$\, is said to converge to
\,$x\,\in\,X$\, if and only if
\[\mathop {\lim }\limits_{n\,\, \to \,\,\infty } \,\mu\,(\,x_{n} \,,\,x \,,\, t\,)\;=\;1 \;\;for \,\,each\,\, t \,>\, 0 \hspace{4.5cm}\]
A sequence $\{\,x_{\,n}\,\}_{\,n}$ in X is called Cauchy sequence if and only if
\[\mathop {\lim }\limits_{n\,\, \to \,\,\infty } \,\mu\,(\,x_{n} \,,\,x_{n\,+\,p} \,,\, t\,)\;=\;1\;for \,\,each\,\, t \,>\, 0 \;and\; p \,=\, 1 \,,\, 2 \,,\, 3 \,,\, \cdots\]
A fuzzy metric space \,$(\,X \,,\, \mu \,,\, \ast\,)$\, is said to be complete if and only if every Cauchy sequence in \,$X$\, is convergent in \,$X$ .
\end{Definition}


\section{Related Fixed Point Theorems}

\begin{Theorem}
Let  \,$(\,X \,,\, \mu \,,\, \ast\,)$\, and \,$(\,Y \,,\, \nu \,,\, \ast\,)$\, be complete fuzzy metric spaces. If \,$T$\, is a continuous mapping of \,$X$\, into \,$Y$\, and \,$S$\, is a mapping of \,$Y$\, into \,$X$\, satisfying the inequalities
\[k\,\mu(\,STx\,,\,STx^{\,'}\,,\,t\,)
\,\geq\, \min\{\,\mu(\,x\,,\,x^{\,'}\,,\,t\,) \,,\,
\mu(\,x\,,\,STx\,,\,t\,) \,,\, \mu(\,x^{\,'}\,,\,STx^{\,'}\,,\,t\,) \,,\,
\nu(\,Tx\,,\,Tx^{\,'}\,,\,t\,)\,\} \]
\[k\,\nu(\,TSy\,,\,TSy^{\,'}\,,\,t\,)
\,\geq\, \min\{\,\nu(\,y\,,\,y^{\,'}\,,\,t\,) \,,\,
\nu(\,y\,,\,TSy\,,\,t\,) \,,\, \nu(\,y^{\,'}\,,\,TSy^{\,'}\,,\,t\,) \,,\,
\mu(\,Sy\,,\,Sy^{\,'}\,,\,t\,)\,\} \]
for all \,$x\,,\,x^{'}$\, in \,$X$\, and \,$y\,,\,y^{'}$\, in \,$Y$\,, where \,$k\,\in\,(\,0\,,\,1\,)$\,, then \,$ST$\, has a unique fixed point \,$z$\, in \,$X$\, and \,$TS$\, has a unique fixed point \,$w$\, in \,$Y$\,. Further, \,$Tz\,=\,w$\, and  \,$Sw\,=\,z$\,
\end{Theorem}

{\bf Proof.} Let \,$x$\, be an arbitrary point in \,$X$\,. Let \,$x_{\,1}\,=\,(ST)^{\,1}x $\,, \,$x_{\,2}\,=\,(ST)^{\,2}x\,=\,(ST)x_{\,1}\,,\cdots\,, x_{\,n}\,=\,(ST)^{\,n}x\,=\,(ST)x_{\,n\,-\,1}$\,
\\and
\\\,$y_{\,1}\,=\,Tx $\,, \,$y_{\,2}\,=\,T(STx)\,=\,Tx_{\,1}\,,\cdots\,, y_{\,n}\,=\,T(ST)^{\,n\,-\,1}x\,=\,Tx_{\,n\,-\,1}$\,
\\for all \,$n\,\in\,N$\,. By induction
\[k\,\mu(\,x_{\,n}\,,\,x_{\,n+1}\,,\,t\,)
\;=\; k\,\mu(\,STx_{\,n\,-\,1}\,,\,STx_{\,n}\,,\,t\,){\hspace{10.5cm}}\]
\[\geq\,\min\{\,\mu(\,x_{\,n\,-\,1}\,,\,x_{\,n}\,,\,t\,) \,,\,
\mu(\,x_{\,n\,-\,1}\,,\,STx_{\,n\,-\,1}\,,\,t\,)
\,,\,\mu(\,x_{\,n}\,,\,STx_{\,n}\,,\,t\,) \,,\, \nu(\,Tx_{\,n\,-\,1}\,,\,Tx_{\,n}\,,\,t\,)\,\} \]
\[=\,\min\{\,\mu(\,x_{\,n\,-\,1}\,,\,x_{\,n}\,,\,t\,)\,,\,
\mu(\,x_{\,n\,-\,1}\,,\,x_{\,n}\,,\,t\,)
\,,\,\mu(\,x_{\,n}\,,\,x_{\,n\,+\,1}\,,\,t\,) \,,\,
\nu(\,y_{\,n}\,,\,y_{\,n\,+\,1}\,,\,t\,)\,\}\]
\begin{equation}
=\,\min\{\,\mu(\,x_{\,n\,-\,1}\,,\,x_{\,n}\,,\,t\,)\,,\,
\mu(\,x_{\,n}\,,\,x_{\,n\,+\,1}\,,\,t\,)
\,,\,\nu(\,y_{\,n}\,,\,y_{\,n\,+\,1}\,,\,t\,)\,\}\hspace{0.5cm} \cdots \hspace{0.7cm}
\end{equation}
Hence
\[\mu(\,x_{\,n}\,,\,x_{\,n\,+\,1}\,,\,t\,){\hspace{11.0cm}}\]
\begin{equation}
\geq\,\frac{1}{k}\,\min\{\,\mu(\,x_{\,n\,-\,1}\,,\,x_{\,n}\,,\,t\,)\,,\,
\mu(\,x_{\,n}\,,\,x_{\,n\,+\,1}\,,\,t\,)
\,,\,\nu(\,y_{\,n}\,,\,y_{\,n\,+\,1}\,,\,t\,)\,\}\hspace{0.4cm} \cdots \hspace{0.4cm}
\end{equation}
By putting \,$(\,2\,)$\, in \,$(\,1\,)$\, we obtain that
\[k\,\mu(\,x_{\,n}\,,\,x_{\,n\,+\,1}\,,\,t\,){\hspace{11.3cm}}\]
\[\geq\,\min\{\,\mu(\,x_{\,n\,-\,1}\,,\,x_{\,n}\,,\,t\,)\,,\,
\nu(\,y_{\,n}\,,\,y_{\,n\,+\,1}\,,\,t\,)\,,\,
\frac{1}{k}\,\mu(\,x_{n\,-\,1}\,,\,x_{\,n}\,,\,t\,)\,,{\hspace{2.5cm}}\]
\[{\hspace{6.6cm}}\frac{1}{k}\,\mu(\,x_{\,n}\,,\,x_{\,n\,+\,1}\,,\,t\,)
\,,\frac{1}{k}\,\nu(\,y_{\,n}\,,\,y_{\,n\,+\,1}\,,\,t\,)\,\}\]
\[=\,\min\{\,\mu(\,x_{\,n\,-\,1}\,,\,x_{\,n}\,,\,t\,)\,,\,
\nu(\,y_{\,n}\,,\,y_{\,n\,+\,1}\,,\,t\,)\,,\,
\frac{1}{k}\,\mu(\,x_{\,n}\,,\,x_{\,n\,+\,1}\,,\,t\,)\,\}{\hspace{3.0cm}}\]
\[\vdots{\hspace{13.0cm}}\]
\[\geq\,\min\{\,\mu(\,x_{\,n\,-\,1}\,,\,x_{\,n}\,,\,t\,)\,,\,
\nu(\,y_{\,n}\,,\,y_{\,n\,+\,1}\,,\,t\,)\,,\,
\frac{1}{k^{\,m}}\,\mu(\,x_{\,n}\,,\,x_{\,n\,+\,1}\,,\,t\,)\,\}{\hspace{2.0cm}}\]
Taking $lim$ as \,$m\longrightarrow\,\infty$\,we have
\begin{equation}
k\,\mu(\,x_{\,n}\,,\,x_{\,n\,+\,1}\,,\,t\,)
\geq\,\min\{\,\mu(\,x_{\,n\,-\,1}\,,\,x_{\,n}\,,\,t\,)\,,\,
\nu(\,y_{\,n}\,,\,y_{\,n\,+\,1}\,,\,t\,)\,\}\hspace{0.4cm} \cdots \hspace{0.4cm}
\end{equation}
Again,
\[k\,\nu(\,y_{\,n}\,,\,y_{\,n\,+\,1}\,,\,t\,)
\;\,=\;\, k\,\nu(\,Tx_{\,n\,-\,1}\,,\,Tx_{\,n}\,,\,t\,) \hspace{7.5cm}\]
\[=\,k\,\nu(\,T(\,STx_{\,n\,-\,2}\,)\,,\,T(\,STx_{\,n\,-\,1}\,)\,,\,t\,)
\;\,=\;\, k\,\nu(\,TSy_{\,n\,-\,1}\,,\,TSy_{\,n}\,,\,t\,){\hspace{7.5cm}}\]
\[\geq\,\min\{\,\nu(\,y_{\,n\,-\,1}\,,\,y_{\,n}\,,\,t\,)\,,\,
\nu(\,y_{\,n\,-\,1}\,,\,TSy_{\,n\,-\,1}\,,\,t\,)
\,,\,\nu(\,y_{\,n}\,,\,TSy_{\,n}\,,\,t\,)\,,{\hspace{2.5cm}}\] \[{\hspace{10.0cm}}\mu(\,Sy_{\,n\,-\,1}\,,\,Sy_{\,n}\,,\,t\,)\,\} \]
\[=\,\min\{\,\nu(\,y_{\,n\,-\,1}\,,\,y_{\,n}\,,\,t\,)\,,\,
\nu(\,y_{\,n\,-\,1}\,,\,y_{\,n}\,,\,t\,)
\,,\,\nu(\,y_{\,n}\,,\,y_{\,n\,+\,1}\,,\,t\,)\,,\,
\mu(\,x_{\,n\,-\,1}\,,\,x_{\,n}\,,\,t\,)\,\}
{\hspace{2.0cm}}\]
\[=\,\min\{\,\nu(\,y_{\,n\,-\,1}\,,\,y_{\,n}\,,\,t\,)\,,\,
\nu(\,y_{\,n}\,,\,y_{\,n\,+\,1}\,,\,t\,)\,,
\,\mu(\,x_{\,n\,-\,1}\,,\,x_{\,n}\,,\,t\,)\,\}{\hspace{3.5cm}}\]
\begin{equation}
\Longrightarrow\; \,k\,\nu(\,y_{\,n}\,,\,y_{\,n\,+\,1}\,,\,t\,)
\geq\,\min\{\,\nu(\,y_{\,n\,-\,1}\,,\,y_{\,n}\,,\,t\,)\,,\,
\mu(\,x_{\,n\,-\,1}\,,\,x_{\,n}\,,\,t\,)\,
\}{\hspace{0.1cm}}\cdots
\end{equation}
From \,$(\,3\,)$\, and \,$(\,4\,)$ , by induction we get
\[\mu(\,x_{\,n}\,,\,x_{\,n+1}\,,\,t\,)\,=\,
\frac{1}{k}\min\{\,\mu(\,x_{\,n\,-\,1}\,,\,x_{\,n}\,,\,t\,)
\,,\,\nu(\,y_{\,n}\,,\,y_{\,n\,+\,1}\,,\,t\,)\,\}{\hspace{6.5cm}}\]
\[\vdots{\hspace{7.0cm}}\]
\[\geq\,\frac{1}{k^{\,n}}\,\min\{\,\mu(\,x\,,\,x_{\,1}\,,\,t\,)
\,,\,\nu(\,y_{\,1}\,,\,y_{\,2}\,,\,t\,)
\,\}{\hspace{0.5cm}}\]
We now verify that \,$ x_{n}$\, is a cauchy sequence. Let \,$t_{\,1}\,=\,\frac{t}{p}$\,.
\[\mu\,(\,x_{\,n} \,,\,x_{\,n\,+\,p} \,,\, t\,){\hspace{11.5cm}}\]
\[\geq\;\mu\,(\,x_{\,n} \,,\,x_{\,n\,+\,1} \,,\,t _{\,1}\,)\,\ast\,\mu\,(\,x_{\,n\,+\,1} \,,\,x_{\,n\,+\,2} \,,\, t_{\,1}\,)\,\ast\,\cdots\,\ast\,\mu\,(\,x_{\,n\,+\,p\,-\,1} \,,\,x_{\,n\,+\,p} \,,\, t_{\,1}\,){\hspace{3.3cm}}\]
\[\geq\,\frac{1}{k^{\,n}}\,\min\{\,\mu(\,x\,,\,x_{\,1}\,,\,t\,)
\,,\,\nu(\,y_{\,1}\,,\,y_{\,2}\,,\,t\,)\,\,\} \;\ast\;\; \cdots{\hspace{7.3cm}}\]
\[{\hspace{5.5cm}}\ast\,\frac{1}{k^{\,n\,+\,p\,-\,1}}\,
\min\{\,\mu(\,x\,,\,x_{\,1}\,,\,t\,)\,,\,\nu(\,y_{\,1}\,,\,y_{\,2}\,,\,t\,)\,\}\]
\[\Longrightarrow \,1\,\geq\,\mathop {\lim }\limits_{n\,\, \to \,\,\infty } \,\mu\,(\,x_{\,n} \,,\,x_{\,n\,+\,p} \,,\,t\,)\,{\hspace{8.0cm}}\]
\[{\hspace{1.1cm}}\geq\,\mathop {\lim }\limits_{n\,\, \to \,\,\infty } \frac{1}{k^{\,n}}\,\min\{\,\mu(\,x\,,\,x_{\,1}\,,\,t\,)
\,,\,\nu(\,y_{\,1}\,,\,y_{\,2}\,,\,t\,)\,\,\}\,>\,1{\hspace{7.6cm}}\]
\[\Longrightarrow\mathop {\lim }\limits_{n\,\, \to \,\,\infty } \,\mu\,(\,x_{\,n} \,,\,x_{\,n\,+\,p} \,,\,t\,)\,=\,1{\hspace{11.0cm}}\]
Hence \,$\{x_{n}\}$\, is a Cauchy sequence with a limit \,$z$\, in X
and similarly , \,$ \{y_{n}\}$\, is a cauchy sequence with a limit \,$w$\, in \,$Y$\,.
\\We have on using the continuity of \,$T$\,
\[\,w\,=\,\mathop {\lim }\limits_{n\,\, \to \,\,\infty } \,y_{\,n}\,=\,\mathop {\lim }\limits_{n\,\, \to \,\,\infty }\,Tx_{\,n}\,=\,Tz\,\]
Further,
\[k\,\mu(\,STz\,,\,x_{\,n}\,,\,t\,)\,=\,
k\,\mu(\,STz\,,\,STx_{\,n\,-\,1}\,,\,t\,){\hspace{8.5cm}}\]
\[\geq\,\min\{\,\mu(\,z\,,\,x_{\,n\,-\,1}\,,\,t\,)\,,\,\mu(\,z\,,\,STz\,,\,t\,)
\,,\,\mu(\,x_{\,n\,-\,1}\,,\,x_{\,n}\,,\,t\,)\,,\,\nu(\,Tz\,,\,y_{\,n}\,,\,t\,)\,\} {\hspace{2.0cm}}\]
and on letting \,$n$\, \,$\longrightarrow \,\,\infty$\, we have
\[k\,\mu(\,STz\,,\,z\,,\,t\,)\geq\,\min\{\,\mu(\,z\,,\,STz\,,\,t\,)\,,\,
\nu(\,Tz\,,\,w\,,\,t\,)\,\}\,=\,\mu(\,z\,,\,STz\,,\,t\,) \]
it follows that \,$STz\,=\,z$\,. Hence we have
\[STz\,=\,Sw\,=\,z\]
Now suppose that \,$ST$\, has a second fixed point\,$z^{\,'}$\,. Then
\[k\,\mu(\,z\,,\,z^{\,'}\,,\,t\,)
\;=\; k\,\mu(\,STz\,,\,STz^{\,'}\,,\,t\,){\hspace{10.5cm}}\]
\[\geq\,\min\{\,\mu(\,z\,,\,z^{\,'}\,,\,t\,)\,,\,\mu(\,z\,,\,STz\,,\,t\,)
\,,\,\mu(\,z^{\,'}\,,\,STz\,,\,t\,)\,,\,\nu(\,Tz\,,\,Tz^{\,'}\,,\,t\,)\}
{\hspace{2.5cm}}\]
\[\Longrightarrow\;\; k\,\mu(\,z\,,\,z^{\,'}\,,\,t\,) \;\geq\; \min\,\{\,\mu(\,z\,,\,z^{\,'}\,,\,t\,) \,,\, 1 \,,\,
\nu(\,Tz\,,\,Tz^{\,'}\,,\,t\,) \,\} \hspace{3.5cm}\]
\[\Longrightarrow\;\; k\,\mu(\,z\,,\,z^{\,'}\,,\,t\,)\,\geq\,\nu(\,Tz\,,\,Tz^{\,'}\,,\,t\,)
{\hspace{10.0cm}}\]
But
\[k\,\nu(\,Tz\,,\,Tz^{\,'}\,,\,t\,) \,=\, \nu(\,TSTz\,,\,TSTz^{\,'}\,,\,t\,){\hspace{7.0cm}}\]
\[\geq\; \min\{\,\nu(\,Tz\,,\,Tz^{\,'}\,,\,t\,)\,,\,
\nu(\,Tz\,,\,TSTz\,,\,t\,)\,,\,
\nu(\,Tz^{\,'}\,,\,TSTz^{\,'}\,,\,t\,) \,,\,
\mu(\,STz\,,\,STz^{\,'}\,,\,t\,)\}\]
\[=\; \min\,\{\, \nu(\,Tz\,,\,Tz^{\,'}\,,\,t\,) \;,\; 1 \;,\; 1 \;,\; \mu(\,z\,,\,z^{\,'}\,,\,t\,) \,\} \hspace{8.5cm}\]
\[\Longrightarrow\;\; \nu(\,Tz\,,\,Tz^{\,'}\,,\,t\,) \;\geq\; \frac{1}{k}\, \mu(\,z\,,\,z^{\,'}\,,\,t\,) \hspace{9.5cm}\]
Hence ,
\[\mu(\,z\,,\,z^{\,'}\,,\,t\,)\,\geq\,\frac{1}{k^{\,2}}\,
\mu(\,z\,,\,z^{\,'}\,,\, t\,)\,\geq\,\cdots\,\geq\,\frac{1}{k^{\,n}}\,
\mu(\,z\,,\,z^{\,'}\,,\, t\,)\,\longrightarrow\,\infty\,\,as \,\, n\longrightarrow\infty{\hspace{1.0cm}}\]
\[\Longrightarrow\,1\,\geq\,\mu(\,z\,,\,z^{\,'}\,,\,t\,)\geq\,\mathop {\lim }\limits_{n\,\, \to \,\,\infty }\frac{1}{k^{\,n}}\,
\mu(\,z\,,\,z^{\,'}\,,\,t\,)\,>\,1 {\hspace{9.5cm}}\]
\[\Longrightarrow \,\mu(\,z\,,\,z^{\,'}\,,\,t\,)\,=\,1{\hspace{11.5cm}}\]
which implies \,$z\,=\,z^{\,'}$\,.
\\Similarly , \,$w$\, is the unique fixed point of \,$TS$\,. This completes the proof of the theorem.

\begin{Theorem}Let \,$(\,X \,,\, \mu \,,\, \ast\,)$\, and \,$(\,Y \,,\, \nu \,,\, \ast\,)$\,be two complete fuzzy metric spaces. Let \,$A\,,\,B$\, be mappings of \,$X$\,into \,$Y$\,and let \,$S\,,\,T$\, be mappings of \,$Y$\, into \,$X$\, satisfying the inequalities
\begin{equation} k\;\mu(\,SAx\,,\,TBx^{\,'}\,,\,t\,) \;\geq\;
\frac{f(\,x\,,\,x^{\,'}\,,\,y\,,\,y^{\,'}\,,\,t\,)}
{h(\,x\,,\,x^{\,'}\,,\,y\,,\,y^{\,'}\,,\,t\,)}\hspace{1.4cm}\cdots
\end{equation}
\begin{equation}k\;\nu(\,BSy\,,\,ATy^{\,'}\,,\,t\,) \;\geq\;
\frac{g(\,x\,,\,x^{\,'}\,,\,y\,,\,y^{\,'}\,,\,t\,)}
{h(\,x\,,\,x^{\,'}\,,\,y\,,\,y^{\,'}\,,\,t\,)}\hspace{1.4cm}\cdots
\end{equation}
for all \,$x\,,\,x^{'}$\, in \,$X$\, and \,$y\,,\,y^{'}$\, in \,$Y$\, for which
\,$f(\,x\,,\,x^{\,'}\,,\,y\,,\,y^{\,'}\,,\,t\,)$\, and \\\,$g(\,x\,,\,x^{\,'}\,,\,y\,,\,y^{\,'}\,,\,t\,)
\,<\,h(\,x\,,\,x^{\,'}\,,\,y\,,\,y^{\,'}\,,\,t\,)\,<\,1$\, where
\[f(\,x\,,\,x^{\,'}\,,\,y\,,\,y^{\,'},\,t\,)\,=\,
\min\{\,\mu(\,x\,,\,x^{\,'},\,t\,)\nu(\,Ax\,,\,Bx^{\,'},\,t\,)\,,\,
\mu(\,x\,,\,x^{\,'},\,t\,)\,\mu(\,Sy\,,\,Ty^{\,'},\,t\,)\,,\]
\[{\hspace{4.5cm}}\mu(\,x\,,\,Ty^{\,'},\,t\,)\,\nu(\,Ax\,,\,ATy^{\,'},\,t\,)\,,\, \mu(\,x^{\,'}\,,\,Sy\,,\,t\,)\,\nu(\,Bx^{\,'}\,,\,BSy\,,\,t\,)\,\}\]
\[g(\,x\,,\,x^{\,'}\,,\,y\,,\,y^{\,'},\,t\,)\,=\,
\min\{\,\nu(\,y\,,\,y^{\,'},\,t\,)\mu(\,Sy\,,\,Ty^{\,'},\,t\,)\,,\,
\nu(\,y\,,\,y^{\,'},\,t\,)\,\nu(\,Ax\,,\,Bx^{\,'},\,t\,)\,,\,\]
\[{\hspace{4.5cm}}\nu(\,y\,,\,Bx^{\,'},\,t\,)\,\mu(\,Sy\,,\,TBx^{\,'},\,t\,)\,,\,
\nu(\,y^{\,'}\,,\,Ax\,,\,t\,)\,\mu(\,Ty^{\,'}\,,\,SAx\,,\,t\,)\,\}\]
\[h(\,x\,,\,x^{\,'}\,,\,y\,,\,y^{\,'}\,,\,t\,)\,=\,\min\{\,\nu(\,Ax\,,\,Bx^{\,'},\,t\,)\,,\,
\mu(\,SAx\,,\,TBx^{\,'},\,t\,)\,,{\hspace{0.6cm}}\]
\[{\hspace{5.8cm}}\,\mu(\,Sy\,,\,Ty^{\,'},\,t\,)\,,\,\nu(\,BSy\,,\,ATy^{\,'},\,t\,)\,\}\]
and \,$0\,<\,k\,<\,1$\,. If one of the mappings \,$A\,,\,B\,,\,S$\, and \,$T$\, is continuous , then \,$SA$\, and
\,$TB$\, have a unique common fixed point \,$z$\,in \,$X$\, and \,$BS$\, and \,$AT$\, have a unique common fixed point \,$w$\,in \,$Y$\,. Further , \,$Az\,=\,Bz\,=\,w$\, and \,$Sw\,=\,Tw\,=\,z$\,.
\end{Theorem}
{\bf Proof.}Let \,$x\,=\,x_{\,0}$\, be an arbitrary point in \,$X$\,. let
\[Ax_{\,0}\,=\,y_{\,1}\,,\,Sy_{\,1}\,=\,x_{\,1}\,,\,Bx_{\,1}\,=,\,y_{\,2}\,,
\,Ty_{\,2}\,=\,x_{\,2}\,\,and\,\,Ax_{\,2}\,=\,y_{\,3}\]
and in general let
\[y_{\,2n\,-\,1}\,=\,Ax_{\,2n\,-\,2}\,,\,x_{\,2n\,-\,1}\,=\,Sy_{\,2n\,-\,1}\,,\,
y_{\,2n}\,=\,Bx_{\,2n\,-\,1}\,\,and \,\,x_{\,2n}\,=\,Ty_{\,2n}\]
for \,$n\,=\,1\,,\,2\,,\cdots$\,
\\We will first of all suppose that for some \,$n$ ,
\[h(\,x_{\,2n}\,,\,x_{\,2n\,-\,1}\,,\,y_{\,2n\,-\,1}\,,\,y_{\,2n}\,,\,t\,)
\,{\hspace{9.0cm}}\]
\[=\,\min\{\,\nu(\,Ax_{\,2n}\,,\,Bx_{\,2n\,-\,1},\,t\,)\,,\,
\mu(\,SAx_{\,2n}\,,\,TBx_{\,2n\,-\,1},\,t\,)\,,\mu(\,Sy_{\,2n\,-\,1}\,,\,Ty_{\,2n},\,t\,)
\,,\,\]
\[{\hspace{9.5cm}}\nu(\,BSy_{\,2n\,-\,1}\,,\,ATy_{\,2n},\,t\,)\,\}\]
\[=\,\min\{\,\nu(\,y_{\,2n\,+\,1}\,,\,y_{\,2n},\,t\,)\,,\,
\mu(\,x_{\,2n\,+\,1}\,,\,x_{\,2n},\,t\,)\,,\mu(\,x_{\,2n\,-\,1}\,,\,x_{\,2n},\,t\,)
\,,\,{\hspace{3.5cm}}\]
\[{\hspace{10.5cm}}\nu(\,y_{\,2n}\,,\,y_{\,2n\,+\,1},\,t\,)\,\}\]
\[=\,1{\hspace{13.5cm}}\]
Then putting\,
$x_{\,2n\,-\,1}\,=\,x_{\,2n}\,=\,x_{\,2n\,+\,1}\,=\,z$ \,and\, $y_{\,2n}\,=\,y_{\,2n\,+\,1}\,=\,w$ ,
\\we see that
\begin{equation}
SAz\,=\,TBz\,=\,z\,,\,ATw\,=\,w\,,\,Az\,=\,Bz\,=\,w\,,\,Tw\,=\,z \hspace{0.5cm} \cdots
\end{equation}
from which it follows that \,$Sw \,=\, z \,,\, BSw \,=\, w$ .
\\Similarly ,
\,$h(\,x_{\,2n}\,,\,x_{\,2n\,+\,1}\,,\,y_{\,2n\,+\,1}\,,\,y_{\,2n}\,,\,t\,)\,=\,1$\,
for some \,$n$\, implies that there exist points \,$z$\, in \,$X$\, and \,$w$\, in \,$Y$\, such that
\[SAz\,=\,TBz\,=\,z\,,\,BSw\,=\,ATw\,=\,w\,,\,Az\,=\,Bz\,=\,w\,,\,Sw\,=\,Tw\,=\,z .\]
We will now suppose that
\[h(\,x_{\,2n}\,,\,x_{\,2n\,-\,1}\,,\,y_{\,2n\,-\,1}\,,\,y_{\,2n}\,,\,t\,)\,<\,1\]
and
\[h(\,x_{\,2n}\,,\,x_{\,2n\,+\,1}\,,\,y_{\,2n\,+\,1}\,,\,y_{\,2n}\,,\,t\,)\,<\,1\]
Applying inequality \,$(5)$ , we get
\[k\,\mu(\,x_{\,2n\,+\,1}\,,\,x_{\,2n}\,,\,t\,)\,\;
=\,\; k\,\mu(\,SAx_{\,2n}\,,\,TBx_{\,2n\,-\,1}\,,\,t\,){\hspace{6.5cm}}\]
\[\geq\,\frac{f(\,x_{\,2n}\,,\,x_{\,2n\,-\,1}\,,
\,y_{\,2n\,-\,1}\,,\,y_{\,2n}\,,\,t\,)}
{h(\,x_{\,2n}\,,\,x_{\,2n\,-\,1}\,,\,y_{\,2n\,-\,1}
\,,\,y_{\,2n}\,,\,t\,)}\]
where
\[f(\,x_{\,2n}\,,\,x_{\,2n\,-\,1}\,,\,y_{\,2n\,-\,1}\,,\,y_{\,2n}\,,\,t\,)
{\hspace{9.0cm}}\]
\[=\,\min\{\,\mu(\,x_{\,2n}\,,\,x_{\,2n\,-\,1},\,t\,)\,
\nu(\,Ax_{\,2n}\,,\,Bx_{\,2n\,-\,1},\,t\,)\,,\,
\mu(\,x_{\,2n}\,,\,x_{\,2n\,-\,1},\,t\,)\,\mu(\,Sy_{\,2n\,-\,1}\,,\,
Ty_{\,2n},\,t\,)\,,\,\]
\[\mu(\,x_{\,2n}\,,\,Ty_{\,2n},\,t\,)
\nu(\,Ax_{\,2n}\,,\,ATy_{\,2n},\,t\,)\,,\,
\mu(\,x_{\,2n\,-\,1}\,,\,Sy_{\,2n\,-\,1}\,,\,t\,)
\nu(\,Bx_{\,2n\,-\,1}\,,\,BSy_{\,2n\,-\,1}\,,\,t\,)\,\}\]
\[=\,\min\{\,\mu(\,x_{\,2n}\,,\,x_{\,2n\,-\,1},\,t\,)\,
\nu(\,y_{\,2n\,+\,1}\,,\,y_{\,2n},\,t\,)\,,\,
\mu(\,x_{\,2n}\,,\,x_{\,2n\,-\,1},\,t\,)\,\mu(\,x_{\,2n\,-\,1}\,,\,
x_{\,2n},\,t\,)\,,\,\]
\[\mu(\,x_{\,2n}\,,\,x_{\,2n},\,t\,)\,
\nu(\,y_{\,2n\,+\,1} \,,\, y_{\,2n\,+\,1},\,t\,) \,,\,
\mu(\,x_{\,2n\,-\,1}\,,\,x_{\,2n\,-\,1}\,,\,t\,)\,
\nu(\,y_{\,2n}\,,\,y_{\,2n}\,,\,t\,)\,\}\]
\[=\,\min\{\,\mu(\,x_{\,2n}\,,\,x_{\,2n\,-\,1},\,t\,)
\nu(\,y_{\,2n\,+\,1}\,,\,y_{\,2n},\,t\,)\,,\,
\mu^{\,2}(\,x_{\,2n}\,,\,x_{\,2n\,-\,1},\,t\,)\,,\,1\,,\,1\,\}{\hspace{2.0cm}}\]
\[=\,\min\{\,\mu(\,x_{\,2n}\,,\,x_{\,2n\,-\,1},\,t\,)
\nu(\,y_{\,2n\,+\,1}\,,\,y_{\,2n},\,t\,)\,,\,
\mu^{\,2}(\,x_{\,2n}\,,\,x_{\,2n\,-\,1},\,t\,)\}{\hspace{2.5cm}}\]
and
\[h(\,x_{\,2n}\,,\,x_{\,2n\,-\,1}\,,\,y_{\,2n\,-\,1}
\,,\,y_{\,2n}\,,\,t\,)\,{\hspace{9.0cm}}\]
\[=\,\min\{\,\nu(\,y_{\,2n\,+\,1}\,,\,y_{\,2n},\,t\,)\,,\,
\mu(\,x_{\,2n\,+\,1}\,,\,x_{\,2n},\,t\,)\,,\mu(\,x_{\,2n\,-\,1}\,,\,x_{\,2n},\,t\,)
\,,\, \nu(\,y_{\,2n}\,,\,y_{\,2n\,+\,1},\,t\,)\,\}\]
then
\[k\,\mu(\,x_{\,2n\,+\,1}\,,\,x_{\,2n}\,,\,t\,)\,{\hspace{10.5cm}}\]
\[\geq\,\frac{\min\{\,\mu(\,x_{\,2n}\,,\,x_{\,2n\,-\,1},\,t\,)
\nu(\,y_{\,2n\,+\,1}\,,\,y_{\,2n},\,t\,)\,,\,
\mu^{\,2}(\,x_{\,2n}\,,\,x_{\,2n\,-\,1},\,t\,)\}}
{\min\{\,\nu(\,y_{\,2n\,+\,1}\,,\,y_{\,2n},\,t\,)\,,\,
\mu(\,x_{\,2n\,+\,1}\,,\,x_{\,2n},\,t\,)\,,\mu(\,x_{\,2n\,-\,1}\,,\,x_{\,2n},\,t\,)
\,\}}{\hspace{3.5cm}}\]
from which it follows that
\begin{equation}
k\,\mu(\,x_{\,2n}\,,\,x_{\,2n\,+\,1}\,,\,t\,)\,\geq\,
\min\{\,\mu(\,x_{\,2n\,-\,1}\,,\,x_{\,2n},\,t\,)\,,\,
\nu(\,y_{\,2n}\,,\,y_{\,2n\,+\,1},\,t\,)\,\} \hspace{0.3cm} \cdots
\end{equation}
Using inequality \,$(\,5\,)$\, again, we get
\[k\,\mu(\,x_{\,2n\,-\,1}\,,\,x_{\,2n}\,,\,t\,)\,{\hspace{10.5cm}}\]
\[=\,k\,\mu(\,SAx_{\,2n\,-\,2}\,,\,TBx_{\,2n\,-\,1}\,,\,t\,){\hspace{9.0cm}}\]
\[\geq\,\frac{f(\,x_{\,2n\,-\,2}\,,\,x_{\,2n\,-\,1}\,,\,y_{\,2n\,-\,1}\,,\,
y_{\,2n\,-\,2}\,,\,t\,)}
{h(\,x_{\,2n\,-\,2}\,,\,x_{\,2n\,-\,1}\,,\,y_{\,2n\,-\,1}\,,\,
y_{\,2n\,-\,2}\,,\,t\,)}{\hspace{9.0cm}}\]
\[=\,\frac{\min\{\,\mu(\,x_{\,2n\,-\,2}\,,\,x_{\,2n\,-\,1},\,t\,)
\nu(\,y_{\,2n\,-\,1}\,,\,y_{\,2n},\,t\,)\,,\,
\mu^{\,2}(\,x_{\,2n\,-\,2}\,,\,x_{\,2n\,-\,1},\,t\,)\}}
{\min\{\,\nu(\,y_{\,2n\,-\,1}\,,\,y_{\,2n},\,t\,)\,,\,
\mu(\,x_{\,2n\,-\,1}\,,\,x_{\,2n},\,t\,)\,,\mu(\,x_{\,2n\,-\,1}\,,\,x_{\,2n\,-\,2},\,t\,)
\,\}}{\hspace{3.5cm}}\]
\begin{equation}
k\,\mu(\,x_{\,2n\,-\,1}\,,\,x_{\,2n}\,,\,t\,)\,\geq\,
\min\{\,\mu(\,x_{\,2n\,-\,2}\,,\,x_{\,2n\,-\,1},\,t\,)\,,\,
\nu(\,y_{\,2n\,-\,1}\,,\,y_{\,2n},\,t\,)\,\}
\end{equation}
Again , on using inequality \,$(\,6\,)$\,
\[k\,\nu(\,y_{\,2n}\,,\,y_{\,2n\,+\,1}\,,\,t\,)
\;=\; k\,\nu(\,BSy_{\,2n\,-\,1}\,,\,ATy_{\,2n}\,,\,t\,) \hspace{7.5cm}\]
\[{\hspace{3.6cm}}\geq\,\frac{g(\,x_{\,2n}\,,\,x_{\,2n\,-\,1}
\,,\,y_{\,2n\,-\,1}\,,\,y_{\,2n}\,,\,t\,)}
{h(\,x_{\,2n}\,,\,x_{\,2n\,-\,1}\,,\,y_{\,2n\,-\,1}
\,,\,y_{\,2n}\,,\,t\,)}{\hspace{9.0cm}}\]
where
\[g(\,x_{\,2n}\,,\,x_{\,2n\,-\,1}\,,\,y_{\,2n\,-\,1}\,,\,y_{\,2n}\,,\,t\,)
{\hspace{10.5cm}}\]
\[=\,\min\{\,\nu(\,y_{\,2n\,-\,1}\,,\,y_{\,2n},\,t\,)\,
\mu(\,Sy_{\,2n\,-\,1}\,,\,Ty_{\,2n},\,t\,)\,,\,
\nu(\,y_{\,2n\,-\,1}\,,\,y_{\,2n},\,t\,)\,
\nu(\,Ax_{\,2n}\,,\,Bx_{\,2n\,-\,1},\,t\,)\,,\,\]
\[\nu(\,y_{\,2n\,-\,1}\,,\,Bx_{\,2n\,-\,1},\,t\,)\,
\mu(\,Sy_{\,2n\,-\,1}\,,\,TBx_{\,2n\,-\,1},\,t\,)\,,\, \nu(\,y_{\,2n}\,,\,Ax_{\,2n}\,,\,t\,)\,
\mu(\,Ty_{\,2n}\,,\,SAx_{\,2n}\,,\,t\,)\,\}\]
\[=\,\min\{\,\nu(\,y_{\,2n\,-\,1}\,,\,y_{\,2n},\,t\,)\,
\mu(\,x_{\,2n\,-\,1}\,,\,x_{\,2n},\,t\,)\,,\,
\nu(\,y_{\,2n\,-\,1}\,,\,y_{\,2n},\,t\,)\,\nu(\,y_{\,2n\,+\,1}\,,\,y_{\,2n},\,t\,)\,,\,\]
\[{\hspace{2.0cm}}\nu(\,y_{\,2n\,-\,1}\,,\,y_{\,2n},\,t\,)\,
\mu(\,x_{\,2n\,-\,1}\,,\,x_{\,2n},\,t\,)\,,\, \nu(\,y_{\,2n}\,,\,y_{\,2n\,+\,1}\,,\,t\,)\,
\mu(\,x_{\,2n}\,,\,x_{\,2n\,+\,1}\,,\,t\,)\,\}\]
\[=\,\min\{\,\nu(\,y_{\,2n\,-\,1}\,,\,y_{\,2n},\,t\,)\,
\mu(\,x_{\,2n\,-\,1}\,,\,x_{\,2n},\,t\,)\,,\,
\nu(\,y_{\,2n\,-\,1}\,,\,y_{\,2n},\,t\,){\hspace{3.5cm}}\]
\[{\hspace{3.5cm}}\nu(\,y_{\,2n\,+\,1} \,,\, y_{\,2n},\,t\,)\,,\,
\nu(\,y_{\,2n}\,,\,y_{\,2n\,+\,1}\,,\,t\,)\,
\mu(\,x_{\,2n}\,,\,x_{\,2n\,+\,1}\,,\,t\,)\,\}\]
We then have either
\[g(\,x_{\,2n}\,,\,x_{\,2n\,-\,1}\,,\,y_{\,2n\,-\,1}\,,\,y_{\,2n}\,,\,t\,) \hspace{7.5cm}\]
\[=\; \nu(\,y_{\,2n\,-\,1}\,,\,y_{\,2n},\,t\,)\,\min\,\{\,
\mu(\,x_{\,2n\,-\,1}\,,\,x_{\,2n},\,t\,) \;,\; \nu(\,y_{\,2n\,+\,1}\,,\,y_{\,2n},\,t\,)\,\}\]
or
\[g(\,x_{\,2n}\,,\,x_{\,2n\,-\,1}\,,\,y_{\,2n\,-\,1}\,,\,y_{\,2n}\,,\,t\,) \hspace{7.5cm}\]
\[=\; \nu(\,y_{\,2n\,+\,1}\,,\,y_{\,2n},\,t\,)\,\min\,\{\,
\mu(\,x_{\,2n}\,,\,x_{\,2n\,+\,1},\,t\,) \;,\; \nu(\,y_{\,2n\,-\,1}\,,\,y_{\,2n},\,t\,)\,\}\]
Further,
\[h(\,x_{\,2n}\,,\,x_{\,2n\,-\,1}\,,\,y_{\,2n\,-\,1}\,,\,y_{\,2n}\,,\,t\,)
{\hspace{9.0cm}}\]
\[=\;\min\,\{\,\nu(\,y_{\,2n\,+\,1}\,,\,y_{\,2n},\,t\,) \,,\,
\mu(\,x_{\,2n\,+\,1}\,,\,x_{\,2n},\,t\,) \,,\, \mu(\,x_{\,2n\,-\,1}\,,\,x_{\,2n},\,t\,)
\,\}{\hspace{3.5cm}}\]
\[=\; \min\,\{\,\nu(\,y_{\,2n\,+\,1}\,,\,y_{\,2n},\,t\,) \,,\,
\mu(\,x_{\,2n\,-\,1}\,,\,x_{\,2n},\,t\,)
\,\}{\hspace{7.5cm}}\]
on using inequality \,$(\,8\,)$ . It follows that either
\[k\,\nu(\,y_{\,2n}\,,\,y_{\,2n\,+\,1}\,,\,t\,) \;\geq\;  \nu(\,y_{\,2n}\,,\,y_{\,2n\,-\,1}\,,\,t\,) \hspace{7.5cm}\]
or
\[k\,\nu(\,y_{\,2n}\,,\,y_{\,2n\,+\,1}\,,\,t\,) \;\geq\;
\min\{\,\mu(\,x_{\,2n\,+\,1}\,,\,x_{\,2n},\,t\,) \,,\,
\nu(\,y_{\,2n\,-\,1}\,,\,y_{\,2n},\,t\,)\,\} \hspace{7.5cm}\]
Thus, we have
\begin{equation}
k\,\nu(\,y_{\,2n}\,,\,y_{\,2n\,+\,1}\,,\,t\,) \;\geq\;
\min\{\,\mu(\,x_{\,2n\,+\,1}\,,\,x_{\,2n},\,t\,) \,,\,
\nu(\,y_{\,2n\,-\,1}\,,\,y_{\,2n},\,t\,)\,\} \;\cdots
\end{equation}
Using inequality \,$(\,6\,)$\, again, we get
\[k\,\nu(\,y_{\,2n}\,,\,y_{\,2n\,-\,1}\,,\,t\,) \;=\;
k\,\nu(\,BSy_{\,2n\,-\,1}\,,\,ATy_{\,2n\,-\,2}\,,\,t\,) \hspace{7.5cm}\]
\[{\hspace{0.5cm}}\geq\,\frac{g(\,x_{\,2n\,-\,2}\,,\,x_{\,2n\,-\,1}\,,\,y_{\,2n\,-\,1}
\,,\,y_{\,2n\,-\,2}\,,\,t\,)}
{h(\,x_{\,2n\,-\,2}\,,\,x_{\,2n\,-\,1}\,,\,y_{\,2n\,-\,1}
\,,\,y_{\,2n\,-\,2}\,,\,t\,)}\]
from which it follows that
\begin{equation}
k\,\nu(\,y_{\,2n}\,,\,y_{\,2n\,-\,1}\,,\,t\,) \;\geq\;
\min\{\,\mu(\,x_{\,2n}\,,\,x_{\,2n\,-\,1},\,t\,) \,,\,
\nu(\,y_{\,2n\,-\,2}\,,\,y_{\,2n\,-\,1},\,t\,)\,\}\; \cdots
\end{equation}
It now follows from inequalities \,$(\,8\,)$\, and \,$(\,10\,)$\, that
\[\mu(\,x_{\,n}\,,\,x_{\,n\,+\,1}\,,\,t\,)
\geq\,\frac{1}{k}\;\min\,\{\mu(\,x_{\,n\,-\,1}\,,\,x_{\,n},\,t\,)\,,\,
\frac{1}{k}\;\mu(\,x_{\,n\,+\,1}\,,\,x_{\,n},\,t\,)\,,\,
\frac{1}{k}\;\nu(\,y_{\,n\,-\,1}\,,\,y_{\,n},\,t\,)\}\]
\[{\hspace{0.4cm}}>\; \frac{1}{k}\;\min\,\{\mu(\,x_{\,n\,-\,1}\,,\,x_{\,n},\,t\,)\,,\,
\nu(\,y_{\,n\,-\,1}\,,\,y_{\,n}\,,\,t\,)\}\]
\[\vdots{\hspace{7.5cm}}\]
\[\geq\; \frac{1}{k^{\,n \,-\, 1}}\,\min\,\{\mu(\,x_{\,1}\,,\,x_{\,2}\,,\,t\,)\,,\,
\nu(\,y_{\,1}\,,\,y_{\,2}\,,\,t\,)\}{\hspace{0.5cm}}\]
Let \,$t_{\,1}\,=\,\frac{t}{p}$ .
\[\mu(\,x_{\,n}\,,\,x_{\,n\,+\,p}\,,\,t\,)\;
\geq \;\mu(\,x_{\,n}\,,\,x_{\,n\,+\,1}\,,\,t_{\,1}\,) \; \ast \; \cdots \; \ast \;
\mu(\,x_{\,n\,+\,p\,-\,1}\,,\,x_{\,n\,+\,p} \,,\, t_{\,1}\,){\hspace{5.5cm}}\]
\[\geq\; \frac{1}{k^{\,n \,-\, 1}}\,\min\,\{\mu(\,x_{\,1}\,,\,x_{\,2}\,,\,t_{\,1}\,)\,,\,
\nu(\,y_{\,1}\,,\,y_{\,2}\,,\,t_{\,1}\,)\}\;\ast\; \cdots \,\;\ast{\hspace{5.5cm}}\]
\[{\hspace{1.5cm}}\frac{1}{k^{\,n\,+\,p\,-\,2}}\,\min\,
\{\mu(\,x_{\,1}\,,\,x_{\,2}\,,\, t_{\,1}\,)\,,\,
\nu(\,y_{\,1}\,,\,y_{\,2}\,,\,t_{\,1}\,)\}\]
which implies that
\[\mathop {\lim }\limits_{n\,\, \to \,\,\infty } \mu(\,x_{\,n}\,,\,x_{\,n\,+\,p}\,,\,t\,) \;>\; 1\,\ast\; \cdots \;\ast\,1 \;=\; 1\]
Similarly ,
\[\mathop {\lim }\limits_{n\,\, \to \,\,\infty } \nu(\,y_{\,n}\,,\,y_{\,n\,+\,p}\,,\,t\,) \;=\; 1{\hspace{3.2cm}}\]
$\Longrightarrow$\,\, $\{x_{\,n}\}$\, is a cauchy sequence in \,$X$\, with a limit \,$z$\, and \,$\{y_{\,n}\}$\, is a cauchy sequence in \,$Y$\, with a limit \,$w$\,.
\\Now suppose that \,$A$\, is continuous . Then
\begin{equation}
w\,=\,\mathop {\lim }\limits_{n\,\, \to \,\,\infty }\,y_{\,2n\,+\,1}
\,=\,\mathop {\lim }\limits_{n\,\, \to \,\,\infty }\,Ax_{\,2n}\,=\,Az \hspace{1.3cm}\cdots
\end{equation}
and
\[\mathop {\lim }\limits_{n\,\, \to \,\,\infty }\,
f(\,z\,,\,x_{\,2n\,-\,1}\,,\,w\,,\,y_{\,2n}\,,\,t\,){\hspace{9.0cm}}\]
\[=\; \mathop {\lim }\limits_{n\,\, \to \,\,\infty }\, \min\{\,\mu(\,z\,,\,x_{\,2n\,-\,1},\,t\,)\nu(\,Az\,,\,Bx_{\,2n\,-\,1},\,t\,)\,,\,
\mu(\,z\,,\,x_{\,2n\,-\,1},\,t\,)\,\mu(\,Sw\,,\,Ty_{\,2n},\,t\,)\,,\,\]
\[{\hspace{1.5cm}}\mu(\,z\,,\,Ty_{\,2n},\,t\,)\,\nu(\,Az\,,\,ATy_{\,2n},\,t\,)\,,\,
\mu(\,x_{\,2n\,-\,1}\,,\,Sw\,,\,t\,)\,\nu(\,Bx_{\,2n\,-\,1}\,,\,BSw\,,\,t\,)\,\}\]
\begin{equation}
=\; \min\{\,\mu(\,z\,,\,Sw,\,t\,)\,\nu(\,w\,,\,BSw,\,t\,) \;,\;
\mu^{\,2}\,(\,z\,,\,Sw,\,t\,)\,\}{\hspace{1.3cm}} \cdots\hspace{0.5cm}
\end{equation}
\[\mathop {\lim }\limits_{n\,\, \to \,\,\infty }\,
g(\,z\,,\,x_{\,2n\,-\,1}\,,\,w\,,\,y_{\,2n}\,,\,t\,){\hspace{9.0cm}}\]
\[=\,\min\{\,\nu(\,w\,,\,y_{\,2n},\,t\,)\,\mu(\,Sw\,,\,Ty_{\,2n},\,t\,) \,,\,
\nu(\,w\,,\,y_{\,2n},\,t\,)\,\nu(\,Az\,,\,Bx_{\,2n\,-\,1},\,t\,) \,,\, \hspace{1.5cm}\]
\[{\hspace{1.5cm}}\nu(\,w\,,\,Bx_{\,2n\,-\,1},\,t\,)\,
\mu(\,Sw\,,\,TBx_{\,2n\,-\,1},\,t\,) \,,\,
\nu(\,y_{\,2n}\,,\,Az\,,\,t\,)\,\mu(\,Ty_{\,2n}\,,\,SAz\,,\,t\,)\,\}\]
\begin{equation}
=\; \min\,\{\,\mu(\,z\,,\,Sw,\,t\,) \,,\, \mu(\,z\,,\,SAz,\,t\,)\,\}{\hspace{2.0cm}} \cdots \hspace{2.5cm}
\end{equation}
\begin{equation}
\mathop {\lim }\limits_{n\,\, \to \,\,\infty }\,h(\,z\,,\,x_{\,2n\,-\,1}\,,\,w\,,\,y_{\,2n}\,,\,t\,)
\;=\; \min\{\,\mu(\,Sw\,,\,z,\,t\,) \,,\,
\nu(\,BSw\,,\,w,\,t\,)\,\}
{\hspace{0.2cm}}\cdots \hspace{2.5cm}
\end{equation}
If
\begin{equation}
\min\,\{\,\mu(\,Sw\,,\,z,\,t\,) \,,\,
\nu(\,BSw\,,\,w,\,t\,)\,\} \;=\; 1 \hspace{2.5cm}\cdots\hspace{2.0cm}
\end{equation}
then
\begin{equation}
Sw \,=\, z \,\,,\,\, BSw \,=\, w \,\,,\,\, Bz \,=\, w.
\hspace{2.5cm}\cdots\hspace{1.5cm}
\end{equation}
If it were possible that
\begin{equation}
\min\{\,\mu(\,Sw\,,\,z,\,t\,) \,,\,
\nu(\,BSw\,,\,w,\,t\,)\,\} \;<\; 1
{\hspace{2.0cm}}\cdots\hspace{1.5cm}
\end{equation}
then we have on using inequality \,$(\,5\,)$\, and equations \,$(\,13\,)$\,and \,$(\,15\,)$
\[\mu(\,Sw\,,\,z\,,\,t\,)\,=\,\mathop {\lim }\limits_{n\,\, \to \,\,\infty }\mu(\,SAz\,,\,TBx_{\,2n\,-\,1}\,,\,t\,) \;\geq\; \frac{1}{k}\,\mathop {\lim }\limits_{n\,\, \to \,\,\infty }\,\frac{f(\,z\,,\,x_{\,2n\,-\,1}\,,\,w\,,\,y_{\,2n}\,,\,t\,)}
{h(\,z\,,\,x_{\,2n\,-\,1}\,,\,w\,,\,y_{\,2n}\,,\,t\,)}{\hspace{9.0cm}}\]
\[=\,\frac{\min\{\,\mu(\,z\,,\,Sw,\,t\,)\,\nu(\,w\,,\,BSw,\,t\,) \,,\,
\mu(\,z\,,\,Sw,\,t\,)\}}{k\,\min\{\,\mu(\,Sw\,,\,z,\,t\,)\,,\,
\nu(\,BSw\,,\,w,\,t\,)\,\}} \;\geq\; \frac{1}{k^{\,2}}\,\;\mu(\,Sw\,,\,z,\,t\,){\hspace{11.5cm}}\]
\[\Longrightarrow\,\mu(\,Sw\,,\,z,\,t\,) \;=\; 1{\hspace{10.5cm}}\]
and so \,$Sw\,=\,z$\,. Further , using inequality \,$(\,6\,)$\, and equations \,$(\,14\,)$\, and \,$(\,15\,)$\,, we have equation \,$(\,17\,)$.
\\To complete the proof, We now prove that \,$Tw \;=\; z$ . Then
\[\mathop {\lim }\limits_{n\,\, \to \,\,\infty }\,
f(\,x_{\,2n}\,,\,z\,,\,w\,,\,w\,,\,t\,){\hspace{11.0cm}}\]
\[=\; \mathop {\lim }\limits_{n\,\, \to \,\,\infty }\,
\min\,\{\,\mu(\,x_{\,2n}\,,\,z,\,t\,)\,\nu(\,Ax_{\,2n}\,,\,Bz,\,t\,) \,,\,
\mu(\,x_{\,2n}\,,\,z,\,t\,)\,\mu(\,Sw\,,\,Tw,\,t\,)\,,\,\]
\[{\hspace{1.0cm}}\mu(\,x_{\,2n}\,,\,Tw,\,t\,)\,\nu(\,Ax_{\,2n}\,,\,ATw,\,t\,)
\,,\, \mu(\,z\,,\,Sw\,,\,t\,)\,\nu(\,Bz\,,\,BSw\,,\,t\,)\,\}\]
\begin{equation}
=\; \min\,\{\,\mu(\,z\,,\,Tw,\,t\,)\,\nu(\,w\,,\,ATw,\,t\,) \,,\,
\nu(\,ATw\,,\,w,\,t\,)\,\}{\hspace{1.0cm}}\cdots\hspace{1.0cm}
\end{equation}
If
\[\mathop {\lim }\limits_{n\,\, \to \,\,\infty }\,
h(\,x_{\,2n}\,,\,z\,,\,w\,,\,w\,,\,t\,)
\;=\; \min\{\,\mu(\,z\,,\,Tw,\,t\,)\,,\,\nu(\,w\,,\,ATw,\,t\,)\,\}
\;=\; 1 \,,\hspace{4.5cm}\]
then obviously \,$Tw \;=\; z$ . So we suppose that
\begin{equation}
\mathop {\lim }\limits_{n\,\, \to \,\,\infty }\,
h(\,x_{\,2n}\,,\,z\,,\,w\,,\,w\,,\,t\,)
\;=\; \min\,\{\,\mu(\,z\,,\,Tw,\,t\,) \,,\, \nu(\,w\,,\,ATw,\,t\,) \,\}
\;<\; 1{\hspace{0.2cm}}\cdots
\end{equation}
Then we have on using inequality \,$(\,5\,)$\, and equations \,$(\,19\,)$\, and \,$(\,20\,)$\,
\[\,\mu(\,z\,,\,Tw\,,\,t\,) \,=\, \mathop {\lim }\limits_{n\,\, \to \,\,\infty }\,\mu(\,SAx_{\,2n}\,,\,TBz\,,\,t\,)
\;\geq\; \frac{1}{k}\;\mathop {\lim }\limits_{n\,\, \to \,\,\infty }\,\frac{f(\,x_{\,2n}\,,\,z\,,\,w\,,\,w\,,\,t\,)}
{h(\,x_{\,2n}\,,\,z\,,\,w\,,\,w\,,\,t\,)}{\hspace{9.0cm}}\]
\[=\; \frac{\min\{\,\mu(\,z\,,\,Tw,\,t\,)\,\nu(\,w\,,\,ATw,\,t\,) \,,\,
\nu(\,ATw\,,\,w,\,t\,)\,\}}{k\,\min\{\,\mu(\,z\,,\,Tw,\,t\,)
\,,\, \nu(\,w\,,\,ATw,\,t\,)\,\}}\,
\;\geq\; \frac{1}{k}\;\mu(\,z\,,\,Tw,\,t\,)\]
\[\Longrightarrow\,\mu(\,z\,,\,Tw,\,t\,) \;=\;1 {\hspace{10.5cm}}\]
We must therefore have \,$Tw\,=\,z$\, and equations \,$(\,7\,)$\, again follow.
\\By the symmetry , the same results again hold if one of the mappings \,$B\,,\,S\,,\,T$\, is continuous , instead of \,$A$\,.
\\To prove the uniqueness , suppose that \,$TB$\, and \,$SA$\, have a second
 common fixed point \,$z^{\,'}$\,. Then , using inequality \,$(\,5\,)$\, , we have
\[\,\mu(\,z\,,\,z^{\,'}\,,\,t\,) \;=\; \mu(\,SAz\,,\,TBz^{\,'}\,,\,t\,)
\;\geq\; \frac{1}{k}\,\frac{f(\,z\,,\,z^{\,'}\,,\,Az\,,\,Bz^{\,'}\,,\,t\,)}
{h(\,z\,,\,z^{\,'}\,,\,Az\,,\,Bz^{\,'}\,,\,t\,)}{\hspace{9.0cm}}\]
\[=\,\frac{\min\{\,\mu(\,z\,,\,z^{\,'},\,t\,)\nu(\,w\,,\,Bz^{\,'},\,t\,)\,,\,
\mu^{\,2}(\,z\,,\,z^{\,'},\,t\,)\,,\,\mu(\,z\,,\,z^{\,'},\,t\,)
\nu(\,w\,,\,Az^{\,'},\,t\,)\,\}}
{k\,\min\{\,\nu(\,w\,,\,Bz^{\,'},\,t\,)\,,\,\mu(\,z\,,\,z^{\,'},\,t\,)
\,,\,\nu(\,w\,,\,Az^{\,'},\,t\,)\,\}}\,
{\hspace{6.5cm}}\]
\[\geq\; \frac{1}{k}\; \mu(\,z\,,\,z^{\,'},\,t\,){\hspace{11.5cm}}\]
\[\Longrightarrow\;\; \mu(\,z\,,\,z^{\,'},\,t\,) \;=\; 1{\hspace{10.5cm}}\]
we can prove similarly that \,$w$\, is the unique common fixed point of \,$BS$\, and \,$AT$\,. This completes the proof of the theorem.

\begin{Corollary}Let \,$A\,,\,B\,,\,S$\, and \,$T$\, be self mappings on the complete fuzzy metric space \,$(\,X \,,\, \mu \,,\, \ast\,)$\, satisfying the inequalities
\begin{equation}
{\hspace{1.0cm}}k\,\mu(\,SAx\,,\,TBy\,,\,t\,) \;\geq\;
\frac{f(\,x\,,\,y\,,\,t\,)}
{h(\,x\,,\,y\,,\,t\,)}{\hspace{1.0cm}}\cdots\hspace{1.5cm}
\end{equation}
\begin{equation}
{\hspace{1.0cm}}k\,\mu(\,BSx\,,\,ATy\,,\,t\,)\,\geq\,
\frac{g(\,x\,,\,y\,,\,t\,)}
{h(\,x\,,\,y\,,\,t\,)}{\hspace{1.0cm}}\cdots\hspace{1.5cm}
\end{equation}
for all \,$x\,,\,y$\, in \,$X$\, for which
\,$f(\,x\,,\,y\,,\,t\,)$\,,\,$g(\,x\,,\,y\,,\,t\,)
\,<\,h(\,x\,,\,y\,,\,t\,)\,<\,1$\, where
\[f(\,x\,,\,y,\,t\,)\,=\,
\min\{\,\mu(\,Sx\,,\,Ty,\,t\,)\mu(\,Ax\,,\,BSx,\,t\,)\,,\,
\mu(\,Sx\,,\,TBy,\,t\,)\mu(\,x\,,\,Sx,\,t\,)\,,\]
\[{\hspace{4.2cm}}\mu(\,x\,,\,y,\,t\,)
\mu(\,SAx\,,\,Ty,\,t\,)\,,\,\mu(\,x\,,\,Ty\,,\,t\,)\mu(\,x\,,\,ATy\,,\,t\,)\,\}\]
\[g(\,x\,,\,y,\,t\,)\,=\,
\min\{\,\mu(\,x\,,\,Sx,\,t\,)\mu(\,x\,,\,y,\,t\,)\,,\,
\mu(\,y\,,\,TBy,\,t\,)\mu(\,y\,,\,Ax,\,t\,)\,,{\hspace{1.5cm}}\]
\[{\hspace{2.5cm}}\,\mu(\,SAx\,,\,Ty,\,t\,)
\mu(\,Ax\,,\,By,\,t\,)\,,\,\mu(\,Ax\,,\,ATy\,,\,t\,)\mu(\,SAx\,,\,Sx\,,\,t\,)\,\}\]
\[h(\,x\,,\,y\,,\,t\,)\,=\,\min\{\,\mu(\,Ax\,,\,BSx,\,t\,)\,,\,
\mu(\,x\,,\,SAx,\,t\,)\,,\,\mu(\,Sx\,,\,TBy,\,t\,)\,,{\hspace{2.5cm}}\]
\[{\hspace{11.0cm}}\,\mu(\,By\,,\,ATy,\,t\,)\,\}\]
and \,$0\,<\,k\,<\,1$\,. If one of the mappings \,$A\,,\,B\,,\,S$\, or \,$T$\, is continuous , then \,$SA$\, and
\,$TB$\, have a unique common fixed point \,$u$\, and \,$BS$\, and \,$AT$\, have a unique common fixed point \,$v$\,. Further , \,$Au\,=\,Bu\,=\,v$\, and \,$Sv\,=\,Tv\,=\,u$\,.
\end{Corollary}


\begin{thebibliography}{0}
\bibitem{Aliouche} A. Aliouche and B. Fisher,
\textit{Fixed point theorems for mappings satisfying implicit relation
on two complete and compact metric spaces}, Applied Mathematics and Mechanics, 27 (9)
(2006), 1217$-$1222.
\bibitem{Veeramani} A. George and P. Veeramani.
\textit{ On Some result in fuzzy metric spaces}, Fuzzy Sets and Systems Vol. 64 $(\,1994\,)$ 395$-$399.

\bibitem{Brian} B. Fisher.
\textit{ Related Fixed Points On Two Metric Spaces}, Mathematics Seminar Notes , Vol. 10 $(\,1982\,)$ 17$-$26.

\bibitem{Schweizer}B. Schweizer , A. Sklar,
\textit{Statistical metric space}, Pacific journal of
 mathematics 10 $(\,1960\,)$ 314$-$334.

\bibitem{Rao} K. P. R. Rao, N. Srinivasa Rao T. Ranga Rao and J. Rajendra Prasad,
\textit{Fixed and related fixed point theorems in sequentially compact fuzzy metric spaces}, Int. Journal of Math. Analysis, Vol. 2, 2008, no. 28, 1353$-$1359 1
\bibitem{Rao1} K. P. R. Rao, Abdelkrim Aliouche and G. Ravi Babu,
\textit{Related Fixed Point Theorems in Fuzzy Metric Spaces}, The Journal of Nonlinear Sciences and its Application, 1 (3) (2008), 194$-$202.

\bibitem{zadeh} L. A. Zadeh
\textit{Fuzzy sets}, Information and control 8 $(\,1965\,)$ 338$-$353.

\bibitem{Telci} M. Telci,
\textit{Fixed points on two complete and compact metric spaces}, Applied Mathematics
and Mechanics, 22 (5) (2001), 564$-$568.

\bibitem{Kramosil}O. Kramosil, J. Michalek   ,
\textit{ Fuzzy metric and statisticalmetric spaces},
 Kybernetica 11 $(\,1975\,)$ 326$-$334.

\bibitem{Namdeo}R. K. Namdeo, S. Jain and Brian Fisher..
\textit{ A Related Fixed Point Theorem For Two Pairs Of Mappings On Two Complete Metric spaces}, Hacettepe Journal of Mathematics and Statistics , Vol. 32 $(\,2003\,)$ 07$-$11.

\bibitem{Samanta}T. K. Samanta  and Iqbal H. Jebril ,
\textit{Finite dimentional intuitionistic fuzzy normed linear
space}, Int. J. Open Problems Compt. Math., Vol 2, No. 4
$(\,2009\,)$ 574$-$591.

\end{thebibliography}
\end{document}